\documentclass{amsproc}%
\usepackage{graphicx}
\usepackage{amscd}
\usepackage{amsmath}
\usepackage{amsfonts}
\usepackage{amssymb}%
\theoremstyle{definition}

\numberwithin{equation}{section}

\begin{document}
\title[LIE PROPERTIES IN ASSOCIATIVE ALGEBRAS]{LIE PROPERTIES IN ASSOCIATIVE ALGEBRAS }
\author{Szilvia Homolya}
\address{\noindent\noindent Institute of Mathematics, University of Miskolc, Miskolc,
Hungary 3515}
\email{mathszil@uni-miskolc.hu}
\author{Jen\H{o} Szigeti}
\address{\noindent\noindent Institute of Mathematics, University of Miskolc, Miskolc,
Hungary 3515}
\email{jeno.szigeti@uni-miskolc.hu}
\author{Leon van Wyk}
\address{Department of Mathematical Sciences, Stellenbosch University, P/Bag X1,
Matieland 7602, Stellenbosch, South Africa }
\email{LvW@sun.ac.za}
\author{Michal Ziembowski}
\address{\noindent\noindent Faculty of Mathematics and Information Science, Technical
University of Warsaw, 00-661 Warsaw, Poland}
\email{m.ziembowski@mini.pw.edu.pl}

\begin{abstract}
Let $K$ be a field, then we exhibit two matrices in the full $n\times n$
matrix algebra~$\mathrm{M}_{n}(K)$ which generate $\mathrm{M}_{n}(K)$ as a Lie
$K$-algebra with the commutator Lie product. We also study Lie centralizers of
a not necessarily commutative unitary algebra and obtain results which we hope
will eventually be a step in the direction of, firstly, proving that a
Lie-nilpotent $K$-subspace (or a sub Lie $K$-algebra) of a finite-dimensional
associative algebra over $K$ of index $k$ (say) generates a Lie-nilpotent
associative subalgebra of much higher nilpotency index, and secondly, in the
light of the sharp upper bound for the maximum ($K$-)dimension of a
Lie-nilpotent $K$-subalgebra of~$M_{n}(K)$ of index $k$ obtained in
\cite{myTAMS}, finding an upper bound for the maximum dimension of a
Lie-nilpotent (of index $k$) sub Lie $K$-algebra of~$\mathrm{M}_{n}(K)$.
Finally, the constructive elementary proof of the Skolem-Noether theorem for
the matrix algebra $\mathrm{M}_{n}(K)$ in \cite{SzivW}, in conjunction with
the well-known characteization of Lie automorphisms of~$\mathrm{M}_{n}(K)$ (if
the characteristic of $K$ is different from 2 and 3) in terms of, amongst
others, automorphisms and anti-automorhisms of $\mathrm{M}_{n}(K)$, leads us
to a unifying approach to constructively describe automorphisms and
anti-automorphisms of $\mathrm{M}_{n}(K)$.

\end{abstract}
\maketitle

\section{Introduction and Motivation}

\bigskip

\noindent Throughout the paper an algebra $R$\ means a not necessarily
commutative unitary algebra over a field $K$ (in most of the results $K$ can
be replaced by a commutative unitary ring satisfying certain mild extra
conditions). The centralizer of an element $a\in R$ is denoted by
$\mathrm{Cen}(a)=\{r\in R:ra=ar\}$, and the centre of $R$ by $\mathrm{Z}%
(R)=\{r\in R:rs=sr$ for all $s\in R\}$. Clearly, $\mathrm{Z}(R)\subseteq
\mathrm{Cen}(a)$ are $K$-subalgebras of $R$.

We start with the following simple observation.

\bigskip

\noindent\textbf{1.1. Proposition.} \textit{If the elements }$a_{1}%
,a_{2},\ldots,a_{t}\in R$\textit{ generate }$R$\textit{ as an associative
algebra, then the intersection of their centralizers is trivial, i.e.,}
\[
\mathrm{Cen}(a_{1})\cap\mathrm{Cen}(a_{2})\cap\cdots\cap\mathrm{Cen}%
(a_{t})=\mathrm{Z}(R).
\]

\bigskip

The full $n\times n$ matrix algebra over $K$ is denoted by $\mathrm{M}_{n}%
(K)$. The standard matrix unit in $\mathrm{M}_{n}(K)$ with $1$ in the $(i,j)$
position and zeros in all other positions is denoted by $E_{i,j}$, and $I_{n}$
denotes the $n\times n$ identity matrix.

The fact that $\mathrm{M}_{n}(K)$ can be generated as a $K$-algebra by the two
matrices $E_{n,1}$ and%
\[
S:=E_{1,2}+E_{2,3}+\cdots+E_{n-1,n},
\]
i.e.,
\begin{equation}
\mathrm{M}_{n}(K)=\langle E_{1,1},S\rangle_{K},
\end{equation}

\bigskip

\noindent played a prominent role in \cite{SzivW}, in which a constructive
elementary proof of the Skolem-Noether theorem (see, e.g., \cite{GSza},
\cite{N} and \cite{Sko}) for the matrix algebra~$\mathrm{M}_{n}(K)$,$\ K$ any
field, was given.
To be precise, given a $K$-automorphism~$\varphi$ of $\mathrm{M}_{n}(K)$, an
invertible matrix $A\in\mathrm{M}_{n}(K)$ yielding the conjugation
\[
\varphi(X)=AXA^{-1}%
\]
\noindent for all $X\in\mathrm{M}_{n}(K)$, was constructed from only the two
$\varphi$-images, $\varphi(E_{n,1})$ and~$\varphi(S)$, of the matrices
$E_{n,1}$ and $S$, respectively, and a nonzero vector $\mathbf{a}$ in the
kernel of the matrix $I_{n}-\bigl(\varphi(S)\bigr)^{n-1}\varphi(E_{n,1}%
)\in\mathrm{M}_{n}(K)$, as follows:%

\begin{equation}
\label{Skolem1}A\negthinspace=\negthinspace\left[  \bigl(\varphi
(S)\bigr)^{n-1}\varphi(E_{n,1})\mathbf{a} \mid\bigl(\varphi(S)\bigr)^{n-2}%
\varphi(E_{n,1})\mathbf{a}\mid\negthinspace\cdots\negthinspace\mid
\varphi(S)\varphi(E_{n,1})\mathbf{a} \mid\varphi(E_{n,1})\mathbf{a} \right]  .
\end{equation}

\bigskip

A Lie automorphism $\psi$ of a $K$-algebra $R$ is a one-to-one $K$-linear map
from $R$ onto itself which preserves the commutator Lie product (also called
the Lie bracket in the literature), i.e.,
\[
\psi([x,y])=[\psi(x),\psi(y)],
\]
equivalently,
\[
\psi(xy-yx)=\psi(x)\psi(y)-\psi(y)\psi(x),
\]
for all $x,y\in R$. We note that $\mathrm{M}_{n}(K)$ with the commutator Lie
product plays an exceptional role in the theory of finite dimensional Lie
algebras. The fundamental Ado-Iwasava theorem (see \cite{Ja}) asserts that
every finite-dimensional Lie $K$-algebra can be embedded into $\mathrm{M}%
_{n}(K)$ for some $n\geq1$.

If $K$ is any field of characteristic different from $2$ and $3$, then (see,
e.g., \cite{Hua}, \cite{M1}, \cite{M2} and \cite{M3}) every Lie automorphism
$\psi$ of $\mathrm{M}_{n}(K)$ can be presented as a sum
\begin{equation}
\psi=\sigma+\tau,
\end{equation}
where $\sigma$ is either an automorphism of $\mathrm{M}_{n}(K)$ (as a
$K$-algebra) or the negative of an anti-automorphism of $\mathrm{M}_{n}(K)$,
and $\tau$ is an additive mapping from $\mathrm{M}_{n}(K)$ to $K$ which maps
commutators into zero. In the light of this significant result we apply
(\ref{Skolem1}) in Section 4, where we present a unifying approach to
constructively describe automorphisms and anti-automorphisms of $\mathrm{M}%
_{n}(K)$.

First we show in Section 2, in a vein similar to [\textbf{14}], that the
matrix $E_{1,1}$ and the permutation matrix $P=S+E_{n,1}$ genenate
$\mathrm{M}_{n}(K)$ as a Lie $K$-algebra.

In Section 3 we study Lie centralizers in a (not necessarily commutative)
unitary algebra $R$. We obtain results which we hope will eventually pave the
way towards, firstly, proving that a Lie-nilpotent $K$-subspace (or a sub Lie
$K$-algebra) of a finite-dimensional associative algebra over $K$ of index $k$
(say) generates a Lie-nilpotent associative subalgebra (of much higher
nilpotency index), and secondly, finding an upper bound (perhaps even a sharp
upper bound) for the maximum dimension of a Lie-nilpotent (of index $k$) sub
Lie $K$-algebra of $\mathrm{M}_{n}(K)$ (see Conjecture 3.8). In this context
the sharp upper bound for the maximum dimension of a Lie-nilpotent
$K$-subalgebra of $\mathrm{M}_{n}(K)$ of index $k\geq1$ is important (see
\cite{myTAMS}).



\section{Two matrices generating $\mathrm{M}_{n}(K)$ as a Lie algebra}

\medskip

\noindent We shall make use of the well known multiplication rule of standard
matrix units:%
\[
E_{i,j}E_{k,l}=\left\{
\begin{array}
[c]{c}%
E_{i,l}\text{ if }j=k;\\
0\text{ if }j\neq k.
\end{array}
\right.
\]
The permutation matrix $P\in\mathrm{M}_{n}(K)$ is defined as follows:%
\[
P=E_{1,2}+E_{2,3}+\cdots+E_{n-1,n}+E_{n,1}.
\]

We now show that $\mathrm{M}_{n}(K)$ can be generated as a Lie $K$-algebra by
two matrices.

\bigskip

\noindent\textbf{2.1. Theorem.}\textit{ The matrices }$P$\textit{ and
}$E_{1,1}$\textit{ generate }$\mathrm{M}_{n}(K)$\textit{ as a Lie algebra with
the commutator Lie product.}

\bigskip

\begin{proof}
Let $\mathcal{G}=\left\langle P,E_{1,1}\right\rangle _{\mathrm{Lie}}$
denote the Lie subalgebra of $\mathrm{M}_{n}(K)$ generated by the matrices $P$
and $E_{1,1}$. Clearly,%
\[
E_{1,2}-E_{n,1}=E_{1,1}P-PE_{1,1}=[E_{1,1},P]\in\mathcal{G}%
\]
and%
\[
E_{1,2}+E_{n,1}=E_{1,1}(E_{1,2}-E_{n,1})-(E_{1,2}-E_{n,1})E_{1,1}%
=[E_{1,1},E_{1,2}-E_{n,1}]\in\mathcal{G}%
\]
ensure that%
\[
E_{1,2}=\frac{1}{2}\Big((E_{1,2}+E_{n,1})+(E_{1,2}-E_{n,1})\Big)\in\mathcal{G}%
\]
and%
\[
E_{n,1}=\frac{1}{2}\Big((E_{1,2}+E_{n,1})-(E_{1,2}-E_{n,1})\Big)\in
\mathcal{G}.
\]

Starting from $E_{1,2}\in\mathcal{G}$, assume that $E_{1,j}\in\mathcal{G}$ for
some $2\leq j\leq n-1$. Using%
\[
S=E_{1,2}+E_{2,3}+\cdots+E_{n-1,n}=P-E_{n,1}\in\mathcal{G},
\]
we obtain that $E_{1,j+1}=[E_{1,j},S]\in\mathcal{G}$. Therefore, it follows
that
\[
E_{1,1},E_{1,2},E_{1,3},\ldots,E_{1,n}\in\mathcal{G}.
\]

Next, starting from $E_{n,1}\in\mathcal{G}$, assume that $E_{i,1}%
\in\mathcal{G}$ for some $3\leq i\leq n$. Now%
\[
E_{i-1,1}-E_{i,2}=SE_{i,1}-E_{i,1}S=[S,E_{i,1}]\in\mathcal{G}%
\]
and%
\[
E_{i,2}=E_{i,1}E_{1,2}-E_{1,2}E_{i,1}=[E_{i,1},E_{1,2}]\in\mathcal{G}%
\]
give that%
\[
E_{i-1,1}=(E_{i-1,1}-E_{i,2})+E_{i,2}\in\mathcal{G}.
\]
Consequently, we have that
\[
E_{n,1},E_{n-1,1},\ldots,E_{2,1},E_{1,1}\in\mathcal{G}.
\]

Finally, if $i\neq j$, then
\[
E_{i,j}=E_{i,1}E_{1,j}-E_{1,j}E_{i,1}=[E_{i,1},E_{1,j}]\in\mathcal{G},
\]
and if $i=j$, then%
\[
E_{i,i}=E_{1,1}+(E_{i,1}E_{1,i}-E_{1,i}E_{i,1})=E_{1,1}+[E_{i,1},E_{1,i}%
]\in\mathcal{G}.
\]

Thus we have that $E_{i,j}\in\mathcal{G}$ for all $1\leq i,j\leq n$, whence we
conclude that $\mathcal{G}=\mathrm{M}_{n}(K)$.
\end{proof}

\bigskip

\noindent\textbf{2.2. Remark.} We note that in the above theorem $K$ can be a
commutative unitary ring such that $\frac{1}{2}\in K$. An other observation is
that the Lie generation of $\mathrm{M}_{n}(K)$\ is much stronger than the
associative generation. Indeed, $E_{1,1}=SE_{2,1}$ implies that $S$ and
$E_{2,1}$\ also generate $\mathrm{M}_{n}(K)$ as an associative $K$-algebra.
Since $S$ and $E_{2,1}$\ have zero traces, it follows that all matrices in
$\left\langle S,E_{2,1}\right\rangle _{\mathrm{Lie}}$ have zero traces and
$\left\langle S,E_{2,1}\right\rangle _{\mathrm{Lie}}\neq\mathrm{M}_{n}(K)$.

\section{The Lie centralizer}

\medskip

\noindent For a sequence $x_{1},x_{2},\ldots,x_{m}$ of elements in a not
necessarily commutative unitary algebra $R$ over a field (or commutative ring)
$K$ with unity we use the notation $[x_{1},x_{2},\ldots,x_{m}]_{m}$\ for the
left normed commutator (or Lie) product:%
\[
\lbrack x_{1}]_{1}=x_{1}\text{ and }[x_{1},x_{2},\ldots,x_{m}]_{m}%
=[\ldots\lbrack\lbrack x_{1},x_{2}],x_{3}],\ldots,x_{m}].
\]

\noindent The $k$-th Lie centralizer of a subset $H\subseteq R$ is%
\[
\mathrm{L}_{k}(H)=\big\{r\in R:[r,x_{1},\ldots,x_{k}]_{k+1}=0\text{ for all
}x_{i}\in H\text{, }1\leq i\leq k\big\},
\]
a $K$-subspace (submodule) of $R$.

As a consequence of $[rs,x_{1}]=[r,sx_{1}]+[s,x_{1}r]$, we can see that the
containment%
\[
\{shr:s,r\in R\text{ and }h\in H\}\subseteq H
\]
implies that $\mathrm{L}_{k}(H)$ is a (unitary) $K$-subalgebra of $R$.
Clearly,%
\[
\underset{h\in H}{\cap}\mathrm{Cen}(h)=\mathrm{L}_{1}(H)\subseteq
\mathrm{L}_{2}(H)\subseteq\cdots\subseteq\mathrm{L}_{k}(H)\subseteq
\mathrm{L}_{k+1}(H)\subseteq\cdots
\]
follows from%
\[
\lbrack r,x_{1},\ldots,x_{k},x_{k+1}]_{k+2}=[[r,x_{1},\ldots,x_{k}%
]_{k+1},x_{k+1}].
\]
The $\omega$-Lie centralizer of $H\subseteq R$ is defined as%
\[
\mathrm{L}_{\omega}(H)=\ \underset{k=1}{\overset{\infty}{{\large \cup}}%
}\mathrm{L}_{k}(H).
\]
A subset $H\subseteq R$ is called Lie-nilpotent of index $k\geq1$ if
$H\subseteq\mathrm{L}_{k}(H)$. A natural further step is the following: $H$ is
called $\omega$-Lie-nilpotent (or almost Lie-nilpotent) if $H\subseteq
\mathrm{L}_{\omega}(H)$.

\bigskip

\noindent\textbf{3.1. Proposition.}\textit{ If }$r\in\mathrm{L}_{k}%
(H)$\textit{ and }$1\leq j\leq k$\textit{, then }%
\[
\lbrack x_{1},\ldots,x_{j},r,x_{j+1},\ldots,x_{k}]_{k+1}=0
\]
\textit{for all }$x_{i}\in H,\ 1\leq i\leq k$\textit{.}

\bigskip

\begin{proof}
It is a well known consequence of the Jacobian identity that in any Lie ring,
$[x_{1},\ldots,x_{j},r]_{j+1}$ can be written as a sum of $2^{j-1}$ terms of
the form
\[
\pm[r, x_{\pi(1)},\ldots,x_{\pi(j)}]_{j+1},
\]
where $\pi$ is some permutation of $\{1,2,\ldots,j\}$. We note that an easy
induction on~$j$ works. It follows that $[x_{1},\ldots,x_{j},r,x_{j+1}%
,\ldots,x_{k}]_{k+1}$ can be written as a sum of~$2^{j-1}$ terms of the form
\[
\pm[r,x_{\pi(1)},\ldots,x_{\pi(j)},x_{j+1},\ldots,x_{k}]_{k+1},
\]
whence $[x_{1},\ldots,x_{j},r,x_{j+1},\ldots,x_{k}]_{k+1}=0$ follows.
\end{proof}

\bigskip

\noindent\textbf{3.2. Proposition.}\textit{ If }$\mathrm{L}_{k}(H)=\mathrm{L}%
_{k+1}(H)$\textit{, then }$\mathrm{L}_{k+1}(H)=\mathrm{L}_{k+2}(H)$\textit{.}

\bigskip

\begin{proof}
For the elements $x_{1}\in H$ and $r_{2}\in\mathrm{L}_{k+2}(H)$ we have
\[
[[r, x_{1}],x_{2},\ldots,x_{k+2}]_{k+3} = [r,x_{1},\ldots,x_{k+2}]_{k+3} = 0
\]
for all $x_{i}\in H, \ 2 \le i \le k+2$. Thus we obtain that $[r,x_{1}]
\in\mathrm{L}_{k+1}(H)$ for all $x_{1}\in H$, whence $[r,x_{1}] \in L_{k}(H)$
and
\[
[r,x_{1},\ldots,x_{k+1}]_{k+2} = [[r,x_{1}],x_{2},\ldots,x_{k},x_{k+1}]_{k+1}
= 0
\]
follow for all $x_{i}\in H, \ 1 \le i \le k + 1$. In view of the above
argument, $r\in\mathrm{L}_{k+1}(H)$ and $\mathrm{L}_{k+2}(H) = \mathrm{L}%
_{k+1}(H)$ can be derived.
\end{proof}

\bigskip

\noindent\textbf{3.3. Proposition.}\textit{ Let }$R$\textit{ be a
finite-dimensional algebra over a field }$K$\textit{ with }$\dim_{K}%
(R)=d$\textit{. Then for any subset }$H\subseteq R$\textit{ we have}
$\mathrm{L}_{\omega}(H)=\mathrm{L}_{d}(H)$\textit{.}

\medskip

\begin{proof}
The finite-dimensionality of $R$ implies that
\[
\{0\}\subseteq\mathrm{L}_{1}(H)\subseteq\mathrm{L}_{2}(H)\subseteq
\cdots\subseteq\mathrm{L}_{k}(H)\subseteq\mathrm{L}_{k+1}(H)\subseteq\cdots
\]
cannot be a strictly ascending infinite chain of $K$-subspaces. In view of
Proposition 3.2, the shape of the above chain is
\[
\{0\}\subset\mathrm{L}_{1}(H)\subset\mathrm{L}_{2}(H)\subset\cdots
\subset\mathrm{L}_{t}(H)=\mathrm{L}_{t+1}(H)=\mathrm{L}_{t+2}(H)=\cdots
\]
for some $t\geq1$ (notice that $1_{R}\in\mathrm{L}_{1}(H)$). Now
\[
t\leq\mathrm{dim}_{K}\big(\mathrm{L}_{t}(H)\big)\leq\mathrm{dim}_{K}(R)=d
\]
and $\mathrm{L}_{\omega}(H)=\mathrm{L}_{d}(H)$ follows.
\end{proof}

\bigskip

\noindent\textbf{3.4. Corollary.}\textit{ Let }$R$\textit{ be a
finite-dimensional algebra over a field }$K$\textit{ with }$\dim_{K}%
(R)=d$\textit{. If }$H\subseteq R$\textit{ is }$\omega$\textit{-Lie-nilpotent
(almost Lie-nilpotent), then }$H$\textit{ is Lie-nilpotent of index }%
$d$\textit{.}

\bigskip

\noindent\textbf{3.5. Theorem.}\textit{ For any subset }$H\subseteq
R$\textit{, we have }$\mathrm{L}_{p}(H)\mathrm{L}_{q}(H)\subseteq
\mathrm{L}_{p+q-1}(H)$\textit{ for all }$p,q\geq1$\textit{, and }%
$\mathrm{L}_{\omega}(H)$\textit{ is a }$K$\textit{-subalgebra of }%
$R$\textit{.}

\bigskip

\begin{proof}
Using an induction on $k\geq1$, we prove that for all $r,s,x_{1},\ldots
,x_{k}\in R$,%
\[
\lbrack rs,x_{1},\ldots,x_{k}]_{k+1}=\negthinspace\negthinspace\negthinspace
\negthinspace\underset{j_{1}<j_{2}<\cdots<j_{k-t}}{\underset{1\leq i_{1}%
<i_{2}<\cdots<i_{t}\leq k}{%
{\displaystyle\sum}
}}\negthinspace\negthinspace\negthinspace\negthinspace\negthinspace
\negthinspace\negthinspace\negthinspace\negthinspace\negthinspace
\negthinspace\negthinspace[r,x_{i_{1}},\ldots,x_{i_{t}}]_{t+1}\cdot\lbrack
s,x_{j_{1}},\ldots,x_{j_{k-t}}]_{k-t+1},\text{ } \text{ }\text{ }
\ast\negthinspace(k)
\]
where the sum is taken over all strictly increasing sequences
\[
1\leq i_{1}<i_{2}<\cdots<i_{t}\leq k \quad\mathrm{and} \quad1\leq j_{1}%
<j_{2}<\cdots<j_{k-t}\leq k,
\]
with $0\leq t\leq k$ and
\[
\{j_{1},j_{2},\ldots,j_{k-t}\}=\{1,2,\ldots,k\}\smallsetminus\{i_{1}%
,,i_{2},\ldots,i_{t}\}.
\]
In the above the empty and the full sequences are allowed with $[r,\varnothing
]_{0+1}= r$ and $[s,\varnothing]_{0+1}=s$.

If $k=1$, then%
\[
\lbrack rs,x_{1}]_{2}=[rs,x_{1}]=r[s,x_{1}]+[r,x_{1}]s=[r,\varnothing
]_{1}\cdot\lbrack s,x_{1}]_{2}+[r,x_{1}]_{2}\cdot\lbrack s,\varnothing]_{1}%
\]
is well known.

Assume that $\ast(k)$ holds for some $k\geq1$. We use
\[
[ab,x_{k+1}]=a[b,x_{k+1}]+[a,x_{k+1}]b
\]
repeatedly in the following calculations:%
\[
\lbrack rs,x_{1},\ldots,x_{k},x_{k+1}]_{k+2}=[[rs,x_{1},\ldots,x_{k}%
]_{k+1},x_{k+1}]
\]%
\[
=\left[  \left(  \underset{j_{1}<j_{2}<\cdots<j_{k-t}}{\underset{1\leq
i_{1}<i_{2}<\cdots<i_{t}\leq k}{%
{\displaystyle\sum}
}}[r,x_{i_{1}},\ldots,x_{i_{t}}]_{t+1}\cdot\lbrack s,x_{j_{1}},\ldots
,x_{j_{k-t}}]_{k-t+1}\right)  ,x_{k+1}\right]
\]
\bigskip%
\[
\negthinspace\negthinspace\negthinspace\negthinspace\negthinspace
\negthinspace\negthinspace\negthinspace\negthinspace\negthinspace
\negthinspace= \underset{j_{1}<j_{2}<\cdots<j_{k-t}}{\underset{1\leq
i_{1}<i_{2}<\cdots<i_{t}\leq k}{%
{\displaystyle\sum}
}}\left[  [r,x_{i_{1}},\ldots,x_{i_{t}}]_{t+1}\cdot\lbrack s,x_{j_{1}}%
,\ldots,x_{j_{k-t}}]_{k-t+1},x_{k+1}\right]
\]
\medskip%
\[
=\left(  \underset{j_{1}<j_{2}<\cdots<j_{k-t}}{\underset{1\leq i_{1}%
<i_{2}<\cdots<i_{t}\leq k}{%
{\displaystyle\sum}
}}[r,x_{i_{1}},\ldots,x_{i_{t}}]_{t+1}\cdot[\lbrack s,x_{j_{1}},\ldots
,x_{j_{k-t}}]_{k-t+1},x_{k+1}] \right)
\]%
\[
+\left(  \underset{j_{1}<j_{2}<\cdots<j_{k-t}}{\underset{1\leq i_{1}%
<i_{2}<\cdots<i_{t}\leq k}{%
{\displaystyle\sum}
}}[[r,x_{i_{1}},\ldots,x_{i_{t}}]_{t+1},x_{k+1}] \cdot\lbrack s,x_{j_{1}%
},\ldots,x_{j_{k-t}}]_{k-t+1}\right)
\]%
\[
=\left(  \underset{j_{1}<j_{2}<\cdots<j_{k-t}}{\underset{1\leq i_{1}%
<i_{2}<\cdots<i_{t}\leq k}{%
{\displaystyle\sum}
}}[r,x_{i_{1}},\ldots,x_{i_{t}}]_{t+1}\cdot\left[  s,x_{j_{1}},\ldots
,x_{j_{k-t}},x_{k+1}\right]  _{k-t+2}\right)
\]%
\[
+\left(  \underset{j_{1}<j_{2}<\cdots<j_{k-t}}{\underset{1\leq i_{1}%
<i_{2}<\cdots<i_{t}\leq k}{%
{\displaystyle\sum}
}}\left[  r,x_{i_{1}},\ldots,x_{i_{t}},x_{k+1}\right]  _{t+2}\cdot\lbrack
s,x_{j_{1}},\ldots,x_{j_{k-t}}]_{k-t+1}\right)
\]
\bigskip%
\[
=\underset{j_{1}^{\prime}<j_{2}^{\prime}<\cdots<j_{(k+1)-m}^{\prime}%
}{\underset{1\leq i_{1}^{\prime}<i_{2}^{\prime}<\cdots<i_{m}^{\prime}\leq
k+1}{%
{\displaystyle\sum}
}}\negthinspace\negthinspace\negthinspace\negthinspace\negthinspace
\negthinspace\negthinspace\negthinspace\negthinspace\negthinspace
\negthinspace\negthinspace\negthinspace\negthinspace\negthinspace
\negthinspace[r,x_{i_{1}^{\prime}},\ldots,x_{i_{m}^{\prime}}]_{m+1}%
\cdot\lbrack s,x_{j_{1}^{\prime}},\ldots,x_{j_{(k+1)-m}^{\prime}}]_{k-m+2}.
\text{ }\text{ }\text{ }\text{ }\text{ }\text{ }\text{ } \ast\negthinspace
(k+1)
\]
\medskip The last equality is a consequence of the fact that a strictly
increasing sequence $1\leq i_{1}^{\prime}<i_{2}^{\prime}<\cdots<i_{m}^{\prime
}\leq k+1$ can appear either as
\[
1\leq i_{1}^{\prime}=i_{1}<i_{2}^{\prime}=i_{2}<\cdots<i_{m}^{\prime}%
=i_{t}\leq k\text{ (with }m=t\text{)}%
\]
or as%
\[
1\leq i_{1}^{\prime}=i_{1}<i_{2}^{\prime}=i_{2}<\cdots<i_{m-1}^{\prime}%
=i_{t}<i_{m}^{\prime}=k+1\text{ (with }m=t+1\text{)}.
\]

If $r\in\mathrm{L}_{p}(H)$, $s\in\mathrm{L}_{q}(H)$, $x_{1},\ldots
,x_{p+q-1}\in H$ and $0\leq t\leq p+q-1$, then either~$p\leq t$ or
$q\leq(p+q-1)-t$, and each summand in%
\[
\lbrack rs,x_{1},\ldots,x_{p+q-1}]_{p+q}
\]%
\[
=\underset{j_{1}<j_{2}<\cdots<j_{k-t}}{\underset{1\leq i_{1}<i_{2}%
<\cdots<i_{t}\leq p+q-1}{%
{\displaystyle\sum}
}}[r,x_{i_{1}},\ldots,x_{i_{t}}]_{t+1}\cdot\lbrack s,x_{j_{1}},\ldots
,x_{j_{(p+q-1)-t}}]_{p+q-t}%
\]
is zero. Indeed, if $p\leq t$, then $r\in\mathrm{L}_{p}(H)$ implies that
$[r,x_{i_{1}},\ldots,x_{i_{t}}]_{t+1}=0$, and if $q\leq(p+q-1)-t$, then
$s\in\mathrm{L}_{q}(H)$ implies that
\[
[s,x_{j_{1}},\ldots,x_{j_{(p+q-1)-t}}]_{p+q-t}=0.
\]
It follows that $rs\in\mathrm{L}_{p+q-1}(H)$.

Since $\mathrm{L}_{p+q-1}(H)\subseteq\mathrm{L}_{\omega}(H)$, we derive that
$\mathrm{L}_{\omega}(H)$ is a $K$-subalgebra of $R$.
\end{proof}

\bigskip

\noindent\textbf{3.6. Remark.} A property $\mathcal{P}$ which is defined for
any finite sequence $x_{1},\ldots,x_{m}$ of elements in $R$ is called
hereditary if $\mathcal{P}$ holds for any subsequence $x_{i_{1}}%
,\ldots,x_{i_{t}}$ with $1\leq i_{1}<i_{2}<\cdots<i_{t}\leq m$. Two typical
examples are $\mathcal{D}$ and $\mathcal{L}$. For a sequence $x_{1}%
,\ldots,x_{m}\in R$ the meaning of $\mathcal{D}$\ is that the elements
$x_{1},\ldots,x_{m}$ are distinct and the meaning of $\mathcal{L}$\ is that
the elements $x_{1},\ldots,x_{m}$ are linearly independent over the base field
$K$.

\noindent The $k$-th Lie centralizer of a subset $H\subseteq R$ with respect
to the property $\mathcal{P}$\ is%
\[
\mathrm{L}_{k}^{\mathcal{P}}(H)=\{r\in R\mid\lbrack r,x_{1},\ldots
,x_{k}]_{k+1}=0\text{ for all }x_{1},\ldots,x_{k}\in H\text{ having property
}\mathcal{P}\}.
\]
Using the same calculations as in the above proof, the following interesting
(and probably far reaching) generalization of Theorem 3.5 can be obtained:
\textit{If }$\mathcal{P}$\textit{\ is a hereditary property, then for any
subset }$H\subseteq R$\textit{, we have }$\mathrm{L}_{p}^{\mathcal{P}%
}(H)\mathrm{L}_{q}^{\mathcal{P}}(H)\subseteq\mathrm{L}_{p+q-1}^{\mathcal{P}%
}(H)$\textit{\ (and the union }$\mathrm{L}_{\omega}^{\mathcal{P}}%
(H)=\cup_{k=1}^{\infty}\mathrm{L}_{k}(H)$\textit{ is a }$K$\textit{-subalgebra
of }$R$\textit{).}

\bigskip

\noindent\textbf{3.7. Remark.} Unfortunately we were not able to prove the following:

\bigskip

\noindent\textit{Let }$R$\textit{\ be a finite dimensional algebra over a
field }$K$\textit{ with }$\dim_{K}(R)=d$\textit{. If }$V\subseteq R$\textit{
is a Lie-nilpotent }$K$\textit{-subspace (or a sub Lie }$K$\textit{-algebra)
of index }$k\geq1$\textit{, then the associative }$K$\textit{-subalgebra
}$\left\langle V\right\rangle _{K}$\textit{ of }$R$\textit{\ generated by }%
$V$\textit{\ is Lie-nilpotent of index }$f(k,d)$\textit{.}

\bigskip

The main result in \cite{myTAMS} states that if $K$ is any field and $R$ is
any Lie-nilpotent $K$-subalgebra of $\mathrm{M}_{n}(K)$ of index $k\geq1$,
then
\[
\mathrm{dim}_{K}(R)\leq g(k+1,n),
\]
where $g(k+1,n)$ is the maximum of
\[
\frac{1}{2}\left(  n^{2}-\sum_{i=1}^{k+1}n_{i}^{2}\right)  +1,
\]
subject to the constraint $\sum_{i=1}^{k+1}n_{i}=n$, with $n_{1}%
,n_{2},...,n_{k+1}$ non-negative integers. To be precise:

\bigskip

\noindent\textbf{Theorem.}\textit{ (see \cite{myTAMS}) If }$R$\textit{ is a
Lie-nilpotent }$K$\textit{-subalgebra of }$\mathrm{M}_{n}(K)$\textit{ of index
}$k\geq1$\textit{, with (according to the Division Algorithm)}
\[
n=(k+1)\left\lfloor \frac{n}{k+1}\right\rfloor +r,\ \ \ 0\leq r<k+1,
\]
\textit{then}%
\[
\frac{1}{2}\Bigg(n^{2}-(k+1-r)\left\lfloor \frac{n}{k+1}\right\rfloor
^{2}-r\bigg(\left\lfloor \frac{n}{k+1}\right\rfloor +1\bigg)^{2}\Bigg)+1
\]
\textit{is a sharp upper bound for }$\mathrm{dim}_{K}(R)$\textit{.}

\bigskip

\noindent Using the above Theorem and the statement formulated in Remark 3.7
would allow to get an upper bound for the maximum dimension of a Lie nilpotent
(of index $k$) sub Lie $K$-algebra of the full matrix algebra $\mathrm{M}%
_{n}(K)$. We had hoped that the foregoing results would lead to a proof of the
following conjecture, but unfortunately we fell short:

\bigskip

\noindent\textbf{3.8. Conjecture:}\textit{ If }$\mathcal{L}\subseteq
\mathrm{M}_{n}(K)$\textit{ is an }$\omega$\textit{-Lie-nilpotent sub Lie }%
$K$\textit{-algebra, then }$\dim_{K}(\mathcal{L})\leq1+\frac{1}{2}(n^{2}%
-n)$\textit{.}

\bigskip

\section{A unifying approach to constructively describe \newline automorphisms
and anti-automorphisms of matrix algebras}

\medskip




The importance of automorphisms and anti-automorphisms of a matrix ring
$\mathrm{M}_{n}(K)$ over a field $K$ is evident.
We apply (\ref{Skolem1}) in this section by presenting a unifying approach to
constructively describe automorphisms and anti-automorphisms of $\mathrm{M}%
_{n}(K)$.

In particuclar, first consider the following setting: for an automorphism $f$
of a field $K$, i.e., for $f\in\mathrm{Aut}(K)$, and for any $X\in
\mathrm{M}_{n}(K)$, let $X_{f}$ denote the matrix obtained from $X$ by
applying $f$ entrywise, i.e., $X_{f}=[x_{i,j}]_{f}=[f(x_{i,j})]$, and let $B$
be any invertible matrix in $\mathrm{M}_{n}(K)$. Then the function
$\beta:\mathrm{M}_{n}(K)\rightarrow\mathrm{M}_{n}(K)$, defined by
\[
\beta(X)=BX_{f}B^{-1}%
\]
for all $X\in\mathrm{M}_{n}(K)$, is a ring automorphism of $\mathrm{M}_{n}%
(K)$, but it need not be a $K$-automorphism of $\mathrm{M}_{n}(K)$. In fact,
it is easily verified that $\beta$ is a $K$-automorphism of $\mathrm{M}%
_{n}(K)$ if and only if $f$ is the identity automorphism of $K$. Nevertheless,
we obtain the following constructive description in the above vein (see also
\cite[Corollary 1.2]{Sem}):

\bigskip

\noindent\textbf{4.1. Proposition.}\textit{ Let }$f\in\mathrm{Aut}(K)$\textit{
(}$K$\textit{ any field), let }$B$\textit{ be any invertible matrix in
}$\mathrm{M}_{n}(K)$\textit{, and let }$\beta:\mathrm{M}_{n}(K)\rightarrow
\mathrm{M}_{n}(K)$\textit{ be the function defined by}%
\[
\beta(X)=BX_{f}B^{-1}%
\]
\noindent\textit{for all }$X\in\mathrm{M}_{n}(K)$\textit{. Then}%
\[
\beta(X)=\overline{B}X_{f}\overline{B}^{-1}%
\]
\noindent\textit{for all }$X\in\mathrm{M}_{n}(K)$\textit{, where }%
$\overline{B}\in\mathrm{M}_{n}(K)$\textit{ is the invertible matrix}%
\[
\overline{B}=\left[  \bigl(\beta(S)\bigr)^{n-1}\beta(E_{n,1})\mathbf{b}%
\mid\bigl(\beta(S)\bigr)^{n-2}\beta(E_{n,1})\mathbf{b}\mid\!\cdots\!\mid
\beta(S)\beta(E_{n,1})\mathbf{b}\mid\beta(E_{n,1})\mathbf{b}\right]  ,
\]
\noindent\textit{and }$\mathbf{b}$\textit{ is a nonzero vector in the kernel
of }$I_{n}-\bigl(\beta(S)\bigr)^{n-1}\beta(E_{n,1})\in\mathrm{M}_{n}%
(K)$\textit{.}

\bigskip

\begin{proof}
Since
\[
\beta(X_{f^{-1}})=BXB^{-1}%
\]
for all $X\in\mathrm{M}_{n}(K)$, it follows that $\alpha:\mathrm{M}%
_{n}(K)\rightarrow\mathrm{M}_{n}(K)$, defined by
\[
\alpha(X)=\beta(X_{f^{-1}})
\]
for all $X\in\mathrm{M}_{n}(K)$, is indeed a $K$-automorphism of
$\mathrm{M}_{n}(K)$. Hence, by (\ref{Skolem1}), we can constructively find an
invertible matrix $\overline{B}$ (say) in $\mathrm{M}_{n}(K)$ such that
\[
\alpha(X)=\overline{B}X\overline{B}^{-1}%
\]
\noindent for all $X\in\mathrm{M}_{n}(K)$, where $\overline{B}\in
\mathrm{M}_{n}(K)$ is the invertible matrix%

\[
\overline{B}\!=\!\left[  \!\bigl(\alpha(S)\bigr)^{n-1}\!\alpha(E_{n,1}%
)\mathbf{b}\mid\bigl(\alpha(S)\bigr)^{n-2}\!\alpha(E_{n,1})\mathbf{b}%
\mid\!\cdots\!\mid\alpha(S)\alpha(E_{n,1})\mathbf{b}\mid\alpha(E_{n,1}%
)\mathbf{b}\right]  ,
\]
\noindent with $\mathbf{b}$ a nonzero vector in the kernel of the matrix
$I_{n}-\bigl(\alpha(S)\bigr)^{n-1}\alpha(E_{n,1})$ in $\mathrm{M}_{n}(K)$.
Since $\alpha(S)=\beta(S_{f^{-1}})$ and $\alpha(E_{n,1})=\beta\bigl((E_{n,1}%
)_{f^{-1}}\bigr)$, and since every entry of $S$ and $E_{n,1}$ is $0$ or $1$,
with $f\in\mathrm{Aut}(K)$, we have that $\alpha(S)=\beta(S)$ and
$\alpha(E_{n,1})=\beta(E_{n,1})$. Therefore,%

\[
\beta(X)=\beta\bigl((X_{f})_{f^{-1}}\bigr)=\alpha(X_{f})=\overline
BX_{f}\overline B^{-1}
\]

\noindent for all $X\in\mathrm{M}_{n}(K)$, where
\[
\overline{B}\!=\!\left[  \!\bigl(\beta(S)\bigr)^{n-1}\beta(E_{n,1}%
)\mathbf{b}\mid\bigl(\beta(S)\bigr)^{n-2}\beta(E_{n,1})\mathbf{b}\mid
\!\cdots\!\mid\beta(S)\beta(E_{n,1})\mathbf{b}\mid\beta(E_{n,1})\mathbf{b}%
\right]  ,
\]
\noindent with $\mathbf{b}$ a nonzero vector in the kernel of $I_{n}%
-\bigl(\beta(S)\bigr)^{n-1}\beta(E_{n,1})\in\mathrm{M}_{n}(K)$.
\end{proof}

\bigskip

For our purposes we state explicitly a result from \cite{Sem}, using our notation:

\bigskip

\noindent\textbf{4.2. Corollary.}\textit{ }(\cite[Corollary 1.2]{Sem})\textit{
Let }$K$\textit{ be an arbitrary field, and let }$\phi:\mathrm{M}%
_{n}(K)\rightarrow\mathrm{M}_{n}(K)$\textit{ be a bijective additive function
satisfying }$\phi(XY)=\phi(X)\phi(Y)$\textit{ for all }$X,Y\in\mathrm{M}%
_{n}(K)$\textit{. Then there exists an automorphism }$f$\textit{ of the field
}$K$\textit{ and an invertible matrix }$A\in\mathrm{M}_{n}(K)$\textit{ such
that}
\[
\phi(X)=AX_{f}A^{-1}%
\]
\noindent\textit{for all }$X\in\mathrm{M}_{n}(K)$\textit{.}

\bigskip

Next, combining Proposition 4.1 and Corollary 4.2, and denoting the transpose
of a matrix $X\in\mathrm{M}_{n}(K)$ by $X^{\top}$, we also obtain the
following constructive and explicit description of an invertible matrix
yielding any anti-automorphism of $\mathrm{M}_{n}(K)$. In this regard it is
noteworthy that, just as $S$ and $E_{n,1}$ generate $\mathrm{M}_{n}(K)$ as a
$K$-algebra, so do their transposes $S^{\top}$ and $E_{1,n}$, respectively.

\bigskip

\noindent\textbf{4.3. Theorem.}\textit{ If }$\phi$\textit{ is a ring
anti-automorphism of }$\mathrm{M}_{n}(K)$\textit{, then}%
\[
\phi(X)=\overline{A}X_{f}^{\top}\overline{A}^{-1}%
\]
\noindent\textit{for all }$X\in\mathrm{M}_{n}(K)$\textit{, where }%
$\overline{A}\in\mathrm{M}_{n}(K)$\textit{ is the invertible matrix}%
\[
\overline{A}\!=\!\left[  \!\bigl(\phi(S^{\top})\bigr)^{n-1}\!\phi
(E_{1,n})\mathbf{a}\!\mid\!\bigl(\phi(S^{\top})\bigr)^{n-2}\!\phi
(E_{1,n})\mathbf{a}\!\mid\!\cdots\!\mid\!\phi(S^{\top})\phi(E_{1,n}%
)\mathbf{a}\!\mid\!\phi(E_{1,n})\mathbf{a}\right]  \!,
\]
\noindent\textit{with }$\mathbf{a}$\textit{ a nonzero vector in the kernel of
}$I_{n}-\bigl(\phi(S^{\top})\bigr)^{n-1}\phi(E_{1,n})\in\mathrm{M}_{n}%
(K)$\textit{.}

\bigskip

\begin{proof}
Let $\mathcal{T}$ denote the transposition map $X\mapsto X^{\top}$ on
$\mathrm{M}_{n}(K)$. Since $\mathcal{T}$ is also a ring anti-automorphism of
$\mathrm{M}_{n}(K)$, the composition $\phi\circ\mathcal{T}$ is a ring
automorphism of $\mathrm{M}_{n}(K)$, and so by Corollary 4.2, there is an
automorphism~$f$ of $K$ and an invertible matrix $A\in\mathrm{M}_{n}(K)$ such
that
\[
(\phi\circ\mathcal{T})(X)=AX_{f}A^{-1}%
\]
\noindent for all $X\in\mathrm{M}_{n}(K)$. Hence, by Proposition 4.1,
\[
(\phi\circ\mathcal{T})(X)=\overline{A}X_{f}\overline{A}^{-1}%
\]
\noindent for all $X\in\mathrm{M}_{n}(K)$, where%

\begin{align*}
\overline A  &  = \Bigl[\bigl((\phi\circ\mathcal{T})(S)\bigr)^{n-1}%
\bigl((\phi\circ\mathcal{T})(E_{n,1})\bigr)\mathbf{a} \mid\bigl((\phi
\circ\mathcal{T})(S)\bigr)^{n-2}\bigl((\phi\circ\mathcal{T})(E_{n,1}%
)\bigr)\mathbf{a}\mid\cdots\\
&  \qquad\qquad\qquad\qquad\qquad\quad\quad\cdots\mid\bigl((\phi
\circ\mathcal{T})(S)\phi(E_{n,1})\bigr)\mathbf{a} \mid\bigl((\phi
\circ\mathcal{T})(E_{n,1})\bigr)\mathbf{a}\Bigr]\\
&  = \Bigl[\bigl(\phi(S^{\top})\bigr)^{n-1}\phi(E_{1,n})\mathbf{a}
\mid\bigl(\phi(S^{\top})\bigr)^{n-2}\phi(E_{1,n})\mathbf{a}\mid\cdots\\
&  \qquad\qquad\qquad\qquad\qquad\qquad\qquad\qquad\qquad\cdots\mid
\phi(S^{\top})\phi(E_{1,n})\mathbf{a} \mid\phi(E_{1,n})\mathbf{a}\Bigr].
\end{align*}
(Here $S^{\top} = E_{2,1} + E_{3,2} + \cdots+ E_{n,n-1}$, and $\mathbf{a}$ is
a nonzero vector in the kernel of $I_{n}-\bigl(\phi(S^{\top})\bigr)^{n-1}%
\phi(E_{1,n})$.)

In particular,
\[
(\phi\circ\mathcal{T})(X^{\top})=AX_{f}^{\top}A^{-1},
\]
i.e.,
\[
\phi(X)=AX_{f}^{\top}A^{-1}%
\]
for all $X\in\mathrm{M}_{n}(K)$.
\end{proof}

\bigskip

We illustrate the construction of $\bar{A}$ in Theorem 4.3 with the
(canonical) symplectic involution as a special case of an anti-automorphism
$\phi$.

\bigskip

\noindent\textbf{Example.} Consider the symplectic involution $\phi$ on
$\mathrm{M}_{8}(K)$ (see, e.g.,~\cite{DHill} or~\cite{RSch}), i.e. $\phi$ is
the anti-automorphism of $\mathrm{M}_{8}(K)$ defined by
\[
\phi\left(  \left[
\begin{array}
[c]{c|c}%
U & P\\\hline
Q & V
\end{array}
\right]  \right)  =\left[
\begin{array}
[c]{r|r}%
V^{\top} & -P^{\top}\\\hline
-Q^{\top} & U^{\top}\\
&
\end{array}
\right]
\]
for all $U,V,P,Q\in\mathrm{M}_{4}(K)$. In order to construct $\overline A$
above, we need certain powers of $\phi(S^{\top})$. Since $S^{\top} = E_{2,1} +
E_{3,2} + \cdots+ E_{8,7}$, we have
\begin{equation}
\label{cumbersome}\phi(S^{\top})= E_{1,2} + E_{2,3} + E_{3,4} + E_{5,6} +
E_{6,7} + E_{7,8} - E_{8,1},
\end{equation}
\noindent i.e.,
\[
S^{\top}=\left[
\begin{array}
[c]{cccc|cccc}%
0 &  &  &  &  &  &  & \\
1 & 0 &  &  &  &  &  & \\
& 1 & 0 &  &  &  &  & \\
&  & 1 & 0 &  &  &  & \\\hline
&  &  & 1 & 0 &  &  & \\
&  &  &  & 1 & 0 &  & \\
&  &  &  &  & 1 & 0 & \\
&  &  &  &  &  & 1 & 0\\
&  &  &  &  &  &  &
\end{array}
\right]  , \quad\phi(S^{\top})=\left[
\begin{array}
[c]{cccc|cccc}%
\negthinspace\negthinspace\negthinspace\negthinspace0 & \negthinspace
\negthinspace\negthinspace\negthinspace1 &  &  &  &  &  & \\
\negthinspace\negthinspace\negthinspace\negthinspace & \negthinspace
\negthinspace\negthinspace\negthinspace0 & 1 &  &  &  &  & \\
\negthinspace\negthinspace\negthinspace\negthinspace & \negthinspace
\negthinspace\negthinspace\negthinspace & 0 & 1 &  &  &  & \\
\negthinspace\negthinspace\negthinspace\negthinspace & \negthinspace
\negthinspace\negthinspace\negthinspace &  & 0 &  &  &  & \\\hline
\negthinspace\negthinspace\negthinspace\negthinspace & \negthinspace
\negthinspace\negthinspace\negthinspace &  &  & 0 & 1 &  & \\
\negthinspace\negthinspace\negthinspace\negthinspace & \negthinspace
\negthinspace\negthinspace\negthinspace &  &  &  & 0 & 1 & \\
\negthinspace\negthinspace\negthinspace\negthinspace & \negthinspace
\negthinspace\negthinspace\negthinspace &  &  &  &  & 0 & 1\\
\negthinspace\negthinspace\negthinspace\negthinspace-1 & \negthinspace
\negthinspace\negthinspace\negthinspace &  &  &  &  &  & 0\\
&  &  &  &  &  &  &
\end{array}
\right]  .
\]

Instead of expressing the higher powers $(S^{\top})^{i}$ and $\bigl(\phi
(S^{\top})\bigr)^{i}, \ i=2,3,\ldots,7$, in the form of an expressions as in
(\ref{cumbersome}), which can obviously be done relatively easily, we have
found the resulting expressions in terms of the $E_{i,j}$'s rather cumbersome
to comprehend, and so, although explicit presentations of these matrices, as
above, take considerably more space, we have opted for the latter, since the
resulting presentations are much more illuminating to the reader. Moreover, it
also makes it absolutely clear that this situation for $\mathrm{M}_{8}(K)$ can
be generalized to $\mathrm{M}_{n}(K)$ for any even number $n$.

Thus, we get the following:%

\[
(S^{\top})^{2}\negthinspace\negthinspace=\negthinspace\negthinspace\left[
\begin{array}
[c]{cccc|cccc}%
0 &  &  &  &  &  &  & \\
0 & 0 &  &  &  &  &  & \\
1 & 0 & 0 &  &  &  &  & \\
& 1 & 0 & 0 &  &  &  & \\\hline
&  & 1 & 0 & 0 &  &  & \\
&  &  & 1 & 0 & 0 &  & \\
&  &  &  & 1 & 0 & 0 & \\
&  &  &  &  & 1 & 0 & 0\\
&  &  &  &  &  &  &
\end{array}
\right]  , \ \bigl(\phi(S^{\top})\bigr)^{2}\negthinspace\negthinspace
=\negthinspace\negthinspace\left[
\begin{array}
[c]{rrrr|rrrr}%
\negthinspace\negthinspace\negthinspace\negthinspace0 & \negthinspace
\negthinspace\negthinspace\negthinspace0 & 1 & 0 &  &  &  & \\
\negthinspace\negthinspace\negthinspace\negthinspace & \negthinspace
\negthinspace\negthinspace\negthinspace0 & 0 & 1 &  &  &  & \\
\negthinspace\negthinspace\negthinspace\negthinspace & \negthinspace
\negthinspace\negthinspace\negthinspace & 0 & 0 &  &  &  & \\
\negthinspace\negthinspace\negthinspace\negthinspace & \negthinspace
\negthinspace\negthinspace\negthinspace &  & 0 &  &  &  & \\\hline
\negthinspace\negthinspace\negthinspace\negthinspace & \negthinspace
\negthinspace\negthinspace\negthinspace &  &  & 0 & 0 & 1 & 0\\
\negthinspace\negthinspace\negthinspace\negthinspace & \negthinspace
\negthinspace\negthinspace\negthinspace &  &  &  & 0 & 0 & 1\\
\negthinspace\negthinspace\negthinspace\negthinspace-1 & \negthinspace
\negthinspace\negthinspace\negthinspace &  &  &  &  & 0 & 0\\
\negthinspace\negthinspace\negthinspace\negthinspace0 & \negthinspace
\negthinspace\negthinspace\negthinspace-1 &  &  &  &  &  & 0\\
&  &  &  &  &  &  &
\end{array}
\right]  ,
\]

\[
(S^{\top})^{3}\negthinspace\negthinspace=\negthinspace\negthinspace\left[
\begin{array}
[c]{rrrr|rrrr}%
0 &  &  &  &  &  &  & \\
0 & 0 &  &  &  &  &  & \\
0 & 0 & 0 &  &  &  &  & \\
1 & 0 & 0 & 0 &  &  &  & \\\hline
0 & 1 & 0 & 0 & 0 &  &  & \\
& 0 & 1 & 0 & 0 & 0 &  & \\
&  & 0 & 1 & 0 & 0 & 0 & \\
&  &  & 0 & 1 & 0 & 0 & 0\\
&  &  &  &  &  &  &
\end{array}
\right]  , \ \bigl(\phi(S^{\top})\bigr)^{3}\negthinspace\negthinspace
=\negthinspace\negthinspace\left[
\begin{array}
[c]{rrrr|rrrrr}%
\negthinspace\negthinspace\negthinspace\negthinspace & \negthinspace
\negthinspace\negthinspace\negthinspace & \negthinspace\negthinspace
\negthinspace\negthinspace & 1 &  &  &  &  & \\
\negthinspace\negthinspace\negthinspace\negthinspace & \negthinspace
\negthinspace\negthinspace\negthinspace & \negthinspace\negthinspace
\negthinspace\negthinspace &  &  &  &  &  & \\
\negthinspace\negthinspace\negthinspace\negthinspace & \negthinspace
\negthinspace\negthinspace\negthinspace & \negthinspace\negthinspace
\negthinspace\negthinspace &  &  &  &  &  & \\
\negthinspace\negthinspace\negthinspace\negthinspace & \negthinspace
\negthinspace\negthinspace\negthinspace & \negthinspace\negthinspace
\negthinspace\negthinspace &  &  &  &  &  & \\\hline
\negthinspace\negthinspace\negthinspace\negthinspace0 & \negthinspace
\negthinspace\negthinspace\negthinspace & \negthinspace\negthinspace
\negthinspace\negthinspace &  &  &  &  &  & 1\\
\negthinspace\negthinspace\negthinspace\negthinspace-1 & \negthinspace
\negthinspace\negthinspace\negthinspace0 & \negthinspace\negthinspace
\negthinspace\negthinspace &  &  &  &  &  & \\
\negthinspace\negthinspace\negthinspace\negthinspace & \negthinspace
\negthinspace\negthinspace\negthinspace-1 & \negthinspace\negthinspace
\negthinspace\negthinspace0 &  &  &  &  &  & \\
\negthinspace\negthinspace\negthinspace\negthinspace & \negthinspace
\negthinspace\negthinspace\negthinspace & \negthinspace\negthinspace
\negthinspace\negthinspace-1 & 0 &  &  &  &  &
\end{array}
\right]  ,
\]

\[
(S^{\top})^{4}\negthinspace\negthinspace=\negthinspace\negthinspace\left[
\begin{array}
[c]{rrrr|rrrrrr}
&  &  &  &  &  &  &  &  & \\
&  &  &  &  &  &  &  &  & \\
&  &  &  &  &  &  &  &  & \\
&  &  &  &  &  &  &  &  & \\\hline
1 &  &  &  &  &  &  &  &  & \\
& 1 &  &  &  &  &  &  &  & \\
&  & 1 &  &  &  &  &  &  & \\
&  &  & 1 &  &  &  &  &  & \\
&  &  &  &  &  &  &  &  &
\end{array}
\right]  , \ \bigl(\phi(S^{\top})\bigr)^{4}\negthinspace\negthinspace
=\negthinspace\negthinspace\left[
\begin{array}
[c]{rrrr|rrrrrr}%
\negthinspace\negthinspace\negthinspace\negthinspace & \negthinspace
\negthinspace\negthinspace\negthinspace & \negthinspace\negthinspace
\negthinspace\negthinspace & \negthinspace\negthinspace &  &  &  &  &  & \\
\negthinspace\negthinspace\negthinspace\negthinspace & \negthinspace
\negthinspace\negthinspace\negthinspace & \negthinspace\negthinspace
\negthinspace\negthinspace & \negthinspace\negthinspace &  &  &  &  &  & \\
\negthinspace\negthinspace\negthinspace\negthinspace & \negthinspace
\negthinspace\negthinspace\negthinspace & \negthinspace\negthinspace
\negthinspace\negthinspace & \negthinspace\negthinspace &  &  &  &  &  & \\
\negthinspace\negthinspace\negthinspace\negthinspace & \negthinspace
\negthinspace\negthinspace\negthinspace & \negthinspace\negthinspace
\negthinspace\negthinspace & \negthinspace\negthinspace &  &  &  &  &  &
\\\hline
\negthinspace\negthinspace\negthinspace\negthinspace-1 & \negthinspace
\negthinspace\negthinspace\negthinspace & \negthinspace\negthinspace
\negthinspace\negthinspace & \negthinspace\negthinspace &  &  &  &  &  & \\
\negthinspace\negthinspace\negthinspace\negthinspace & \negthinspace
\negthinspace\negthinspace\negthinspace-1 & \negthinspace\negthinspace
\negthinspace\negthinspace & \negthinspace\negthinspace &  &  &  &  &  & \\
\negthinspace\negthinspace\negthinspace\negthinspace & \negthinspace
\negthinspace\negthinspace\negthinspace & \negthinspace\negthinspace
\negthinspace\negthinspace-1\negthinspace\negthinspace &  &  &  &  &  &  & \\
\negthinspace\negthinspace\negthinspace\negthinspace & \negthinspace
\negthinspace\negthinspace\negthinspace & \negthinspace\negthinspace
\negthinspace\negthinspace & \negthinspace\negthinspace-1 &  &  &  &  &  &
\end{array}
\right]  ,
\]

\[
(S^{\top})^{5}\negthinspace\negthinspace=\negthinspace\negthinspace\left[
\begin{array}
[c]{rrrr|rrrrrr}
&  &  &  &  &  &  &  &  & \\
&  &  &  &  &  &  &  &  & \\
&  &  &  &  &  &  &  &  & \\
&  &  &  &  &  &  &  &  & \\\hline
0 &  &  &  &  &  &  &  &  & \\
1 & 0 &  &  &  &  &  &  &  & \\
& 1 & 0 &  &  &  &  &  &  & \\
&  & 1 & 0 &  &  &  &  &  & \\
&  &  &  &  &  &  &  &  &
\end{array}
\right]  , \ \bigl(\phi(S^{\top})\bigr)^{5}\negthinspace\negthinspace
=\negthinspace\negthinspace\left[
\begin{array}
[c]{rrrr|rrrrrr}%
\negthinspace\negthinspace & \negthinspace\negthinspace\negthinspace
\negthinspace & \negthinspace\negthinspace & \negthinspace\negthinspace &  &
&  &  &  & \\
\negthinspace\negthinspace & \negthinspace\negthinspace\negthinspace
\negthinspace & \negthinspace\negthinspace & \negthinspace\negthinspace &  &
&  &  &  & \\
\negthinspace\negthinspace & \negthinspace\negthinspace\negthinspace
\negthinspace & \negthinspace\negthinspace & \negthinspace\negthinspace &  &
&  &  &  & \\
\negthinspace\negthinspace & \negthinspace\negthinspace\negthinspace
\negthinspace & \negthinspace\negthinspace & \negthinspace\negthinspace &  &
&  &  &  & \\\hline
\negthinspace\negthinspace0 & \negthinspace\negthinspace\negthinspace
\negthinspace-1 & \negthinspace\negthinspace & \negthinspace\negthinspace &  &
&  &  &  & \\
\negthinspace\negthinspace & \negthinspace\negthinspace\negthinspace
\negthinspace0 & \negthinspace\negthinspace-1 & \negthinspace\negthinspace &
&  &  &  &  & \\
\negthinspace\negthinspace & \negthinspace\negthinspace\negthinspace
\negthinspace & \negthinspace\negthinspace0 & \negthinspace\negthinspace-1 &
&  &  &  &  & \\
\negthinspace\negthinspace & \negthinspace\negthinspace\negthinspace
\negthinspace & \negthinspace\negthinspace & \negthinspace\negthinspace0 &  &
&  &  &  & \\
&  &  &  &  &  &  &  &  &
\end{array}
\right]  ,
\]

\[
(S^{\top})^{6}\negthinspace\negthinspace=\negthinspace\negthinspace\left[
\begin{array}
[c]{rrrr|rrrrrr}
&  &  &  &  &  &  &  &  & \\
&  &  &  &  &  &  &  &  & \\
&  &  &  &  &  &  &  &  & \\
&  &  &  &  &  &  &  &  & \\\hline
0 &  &  &  &  &  &  &  &  & \\
0 & 0 &  &  &  &  &  &  &  & \\
1 & 0 & 0 &  &  &  &  &  &  & \\
0 & 1 & 0 & 0 &  &  &  &  &  & \\
&  &  &  &  &  &  &  &  &
\end{array}
\right]  , \ \bigl(\phi(S^{\top})\bigr)^{6}\negthinspace\negthinspace
=\negthinspace\negthinspace\left[
\begin{array}
[c]{rrrr|rrrrrr}%
\negthinspace\negthinspace & \negthinspace\negthinspace\negthinspace
\negthinspace & \negthinspace\negthinspace & \negthinspace\negthinspace &  &
&  &  &  & \\
\negthinspace\negthinspace & \negthinspace\negthinspace\negthinspace
\negthinspace & \negthinspace\negthinspace & \negthinspace\negthinspace &  &
&  &  &  & \\
\negthinspace\negthinspace & \negthinspace\negthinspace\negthinspace
\negthinspace & \negthinspace\negthinspace & \negthinspace\negthinspace &  &
&  &  &  & \\
\negthinspace\negthinspace & \negthinspace\negthinspace\negthinspace
\negthinspace & \negthinspace\negthinspace & \negthinspace\negthinspace &  &
&  &  &  & \\\hline
\negthinspace\negthinspace0 & 0 & \negthinspace\negthinspace\negthinspace
\negthinspace-1 & \negthinspace\negthinspace & \negthinspace\negthinspace &  &
&  &  & \\
\negthinspace\negthinspace & 0 & \negthinspace\negthinspace\negthinspace
\negthinspace0 & \negthinspace\negthinspace-1 & \negthinspace\negthinspace &
&  &  &  & \\
\negthinspace\negthinspace &  & \negthinspace\negthinspace\negthinspace
\negthinspace0 & \negthinspace\negthinspace0 & \negthinspace\negthinspace &  &
&  &  & \\
\negthinspace\negthinspace & \negthinspace\negthinspace\negthinspace
\negthinspace & \negthinspace\negthinspace & \negthinspace\negthinspace0 &  &
&  &  &  & \\
&  &  &  &  &  &  &  &  &
\end{array}
\right]  ,
\]
\vskip0.2truecm \noindent and%

\[
(S^{\top})^{7}\negthinspace\negthinspace=\negthinspace\negthinspace\left[
\begin{array}
[c]{rrrrr|rrrrrr}
&  &  &  &  &  &  &  &  &  & \\
&  &  &  &  &  &  &  &  &  & \\
&  &  &  &  &  &  &  &  &  & \\
&  &  &  &  &  &  &  &  &  & \\\hline
&  &  &  &  &  &  &  &  &  & \\
&  &  &  &  &  &  &  &  &  & \\
&  &  &  &  &  &  &  &  &  & \\
1 &  &  &  &  &  &  &  &  &  & \\
&  &  &  &  &  &  &  &  &  &
\end{array}
\right]  , \ \bigl(\phi(S^{\top})\bigr)^{7}\negthinspace\negthinspace
=\negthinspace\negthinspace\left[
\begin{array}
[c]{rrrr|rrrrr}
&  &  & \negthinspace\negthinspace &  &  &  &  & \\
&  &  & \negthinspace\negthinspace &  &  &  &  & \\
&  &  & \negthinspace\negthinspace &  &  &  &  & \\
&  &  & \negthinspace\negthinspace &  &  &  &  & \\\hline
&  &  & \negthinspace\negthinspace-1 & \negthinspace\negthinspace &  &  &  &
\\
&  &  & \negthinspace\negthinspace & \negthinspace\negthinspace &  &  &  & \\
&  &  & \negthinspace\negthinspace & \negthinspace\negthinspace &  &  &  & \\
&  &  & \negthinspace\negthinspace & \negthinspace\negthinspace &  &  &  & \\
&  &  &  &  &  &  &  &
\end{array}
\right]  .
\]

\noindent Hence, $\phi(S^{\top})\bigr)^{7}\phi(E_{1,8})=(-E_{5,4}%
)(-E_{4,5})=E_{5,5}$, and so
\[
\mathbf{a}:=e_{5}=\left[
\begin{array}
[c]{c}%
0\\
0\\
0\\
0\\
1\\
0\\
0\\
0\\
\end{array}
\right]
\]
\noindent is a nonzero vector in the kernel of
\[
I_{8}-\phi(S^{\top})\bigr)^{7}\phi(E_{1,8})=E_{1,1}+E_{2,2}+E_{3,3}%
+E_{4,4}+E_{6,6}+E_{7,7}+E_{8,8}.
\]
(Here $e_{j}$ denotes the $8\times1$ column vector with $1$ in position $j$,
and $0$ elsewhere.) Therefore, since%

\[
\phi(E_{1,8})\mathbf{a}=-E_{4,5}e_{5}=-e_{4}=\left[
\begin{array}
[c]{c}%
0\\
0\\
0\\
-1\\
0\\
0\\
0\\
0\\
\end{array}
\right]  ,
\]
the foregoing presentations of $\bigl(\phi(S^{\top})\bigr)^{i},\ i=2,3,\ldots
,7$, together with the construction of $\overline{A}$ in Proposition
\ref{proposition3}, yields
\begin{align*}
\overline{A}  &  =\left[  \bigl(\phi(S^{\top})\bigr)^{7}\phi(E_{1,8}%
)\mathbf{a}\mid\bigl(\phi(S^{\top})\bigr)^{6}\phi(E_{1,8})\mathbf{a}\mid
\cdots\mid\phi(S^{\top})\phi(E_{1,n})\mathbf{a}\mid\phi(E_{1,n})\mathbf{a}%
\right] \\
&  =\left[
\begin{array}
[c]{rrrr|rrrr}
&  &  &  & \!\!\!-1 & \!\!\! & \!\!\! & \!\!\!\\
&  &  &  & \!\!\! & \!\!\!-1 & \!\!\! & \!\!\!\\
&  &  &  & \!\!\! & \!\!\!\! & \!\!\!\!-1 & \!\!\!\!\\
&  &  &  & \!\!\! & \!\!\!\! & \!\!\!\! & \!\!\!\!-1\\\hline
1 &  &  &  & \!\!\! & \!\!\!\! & \!\!\!\! & \!\!\!\!\\
& 1 &  &  & \!\!\! & \!\!\!\! & \!\!\!\! & \!\!\!\!\\
&  & 1 &  & \!\!\! & \!\!\!\! & \!\!\!\! & \!\!\!\!\\
&  &  & 1 & \!\!\! & \!\!\!\! & \!\!\!\! & \!\!\!\!\\
&  &  &  &  &  &  &
\end{array}
\right] \\
&  =\left[
\begin{array}
[c]{r|r}%
0 & -I_{4}\\\hline
I_{4} & 0
\end{array}
\right]  ,
\end{align*}
\noindent the latter being the negative of the matrix $y$ on the last line of
the first page of~\cite{RSch}. (Of course, $\overline{A}X^{\top}\overline
{A}^{-1}=(\lambda\overline{A})X^{\top}(\lambda\overline{A})^{-1}$ for every
$0\neq\lambda\in F$.) This concludes the example.

\bigskip

Next, consider the setting following (\ref{Skolem1}), with the only difference
that
the function $\beta:\mathrm{M}_{n}(K)\rightarrow\mathrm{M}_{n}(K)$ is defined
by
\[
\beta(X)=BX_{f}^{\top}B^{-1}%
\]
for all $X\in\mathrm{M}_{n}(K)$ (instead of $\beta(X)=BX_{f}B^{-1}$). Then, as
before, $\beta$ is a ring anti-automorphism of $\mathrm{M}_{n}(K)$, but it
need not be a $K$-anti-automorphism of $\mathrm{M}_{n}(K)$. In this case we
have the following result:

\bigskip

\noindent\textbf{4.4. Corollary.}\textit{ Let }$f\in\mathrm{Aut}(K)$\textit{
(}$K$\textit{ any field), let }$B$\textit{ be any invertible matrix in
}$\mathrm{M}_{n}(K)$\textit{, and let }$\beta:\mathrm{M}_{n}(K)\rightarrow
\mathrm{M}_{n}(K)$\textit{ be the function defined by}
\[
\beta(X)=BX_{f}^{\top}B^{-1}%
\]
\noindent\textit{for all }$X\in\mathrm{M}_{n}(K)$\textit{. Then}%
\[
\beta(X)=\overline{B}X_{f}^{\top}\overline{B}^{-1}%
\]
\noindent\textit{for all }$X\in\mathrm{M}_{n}(K)$\textit{, where }%
$\overline{B}\in\mathrm{M}_{n}(K)$\textit{ is the invertible matrix}%
\[
\overline{B}\!=\!\left[  \!\bigl(\beta(S^{\top})\bigr)^{n-1}\!\beta
(E_{1,n})\mathbf{b}\!\mid\!\bigl(\beta(S^{\top})\bigr)^{n-2}\!\beta
(E_{1,n})\mathbf{b}\!\mid\!\cdots\!\mid\!\beta(S^{\top})\beta(E_{1,n}%
)\mathbf{b}\!\mid\!\beta(E_{1,n})\mathbf{b}\right]
\]
\noindent\textit{and }$\mathbf{b}$\textit{ is a nonzero vector in the kernel
of }$I_{n}-\bigl(\beta(S^{\top})\bigr)^{n-1}\beta(E_{1,n})\in\mathrm{M}%
_{n}(K)$\textit{.}

\bigskip

\begin{proof}
Since
\[
\beta(X_{f^{-1}})=BX^{\top}B^{-1}%
\]
for all $X\in\mathrm{M}_{n}(K)$, it follows that $\alpha:\mathrm{M}%
_{n}(K)\rightarrow\mathrm{M}_{n}(K)$, defined by
\[
\alpha(X)=\beta(X_{f^{-1}})
\]
for all $X\in\mathrm{M}_{n}(K)$, is a $K$-anti-automorphism of $\mathrm{M}%
_{n}(K)$. Hence, by Theorem 4.3, we can constructively find an invertible
matrix $\overline{B}$ (say) in $\mathrm{M}_{n}(K)$ such that
\[
\alpha(X)=\overline{B}X^{\top}\overline{B}^{-1}%
\]
\noindent for all $X\in\mathrm{M}_{n}(K)$, where $\overline{B}\in
\mathrm{M}_{n}(K)$ is the invertible matrix%

\[
\overline B\negthinspace=\negthinspace\negthinspace\left[  \negthinspace
\bigl(\alpha(S^{\top})\bigr)^{n-1}\negthinspace\alpha(E_{1,n})\mathbf{b}%
\negthinspace\mid\negthinspace\bigl(\alpha(S^{\top})\bigr)^{n-2}%
\negthinspace\alpha(E_{1,n})\mathbf{b}\negthinspace\mid\negthinspace
\cdots\negthinspace\mid\negthinspace\alpha(S^{\top})\alpha(E_{1,n}%
)\mathbf{b}\negthinspace\mid\negthinspace\alpha(E_{1,n})\mathbf{b} \right]
\]
\noindent and $\mathbf{b}$ is a nonzero vector in the kernel of the matrix
$I_{n}-\bigl(\alpha(S^{\top})\bigr)^{n-1}\alpha(E_{1,n})$ in~$\mathrm{M}%
_{n}(K)$. We have $\alpha(S^{\top})=\beta(S^{\top}_{f^{-1}})$ and
$\alpha(E_{1,n})=\beta\bigl((E_{n,1})_{f^{-1}}\bigr)$, and so, since every
entry of both $S^{\top}$ and $E_{1,n}$ is $0$ or $1$, and since $f^{-1}%
\in\mathrm{Aut}(K)$, we have that $\alpha(S^{\top})=\beta(S^{\top})$ and
$\alpha(E_{1,n})=\beta(E_{1,n})$. Therefore,%

\[
\beta(X)=\beta\bigl((X_{f})_{f^{-1}}\bigr)=\alpha(X_{f})=\overline BX^{\top
}_{f}\overline B^{-1}
\]

\noindent for all $X\in\mathrm{M}_{n}(K)$, where
\[
\overline{B}\!=\!\!\left[  \!\bigl(\beta(S^{\top})\bigr)^{n-1}\!\beta
(E_{1,n})\mathbf{b}\!\mid\!\bigl(\beta(S^{\top})\bigr)^{n-2}\!\beta
(E_{1,n})\mathbf{b}\!\mid\!\cdots\!\mid\!\beta(S^{\top})\beta(E_{1,n}%
)\mathbf{b}\!\mid\!\beta(E_{1,n})\mathbf{b}\right]  \!\!,
\]
\noindent with $\mathbf{b}$ a nonzero vector in the kernel of $I_{n}%
-\bigl(\beta(S^{\top})\bigr)^{n-1}\beta(E_{1,n})\in\mathrm{M}_{n}(K)$.
\end{proof}

\bigskip


Consider again (1.3). By Proposition 4.1, Corollary 4.2 and Theorem 4.3 we
have an exact description of $\sigma$ in (1.3) in the terms of the images of
generators of $\mathrm{M}_{n}(K)$. Regarding $\tau$, recall that it is an
additive mapping from $\mathrm{M}_{n}(K)$ to $K$ which maps commutators into
zero. With $\mathrm{tr}(X)$ denoting the trace of a matrix $X$ in
$\mathrm{M}_{n}(K)$, we have the following:


\bigskip

\noindent\textbf{4.5. Proposition.}\textit{ Let }$\tau$\textit{ be as in
(1.3), and let }$X=\sum_{i,j=1}^{n}k_{ij}E_{i,j}\in\mathrm{M}_{n}(K)$\textit{.
Then }$\tau(X)=\tau\big(\mathrm{tr}(X)\cdot E_{1,1}\big).$

\bigskip

\begin{proof}
If $i\neq j$, then $[k_{ij}E_{i,j},E_{j,j}]=k_{ij}E_{i,j}$, and so, since
$\tau$ maps commutators to zero, we have $\tau(k_{ij}E_{i,j})=0.$ Hence,
$\tau(X)=\tau\big(\sum_{i=1}^{n}k_{ii}E_{i,i}\big)$. Note also that
$k_{ii}E_{i,i}=[k_{ii}E_{i,1},E_{1,i}]+k_{ii}E_{1,1}$ for every $i$, and so
$\tau(k_{ii}E_{i,i})=\tau(k_{ii}E_{1,1})$. Consequently,
\[
\tau(X)=\tau\Biggl(\biggl(\sum_{i=1}^{n}k_{ii}\biggr)E_{1,1}\Biggr)=\tau
\big(\mathrm{tr}(X)\cdot E_{1,1}\big).
\]

\end{proof}

\bigskip


Unfortunately, we do not seem to be able to describe $\tau\big(\mathrm{tr}%
(X)\cdot E_{1,1}\big)$ any better. In general, if $\psi$ in (1.3) is not a Lie
$K$-automorphism, then we may not have $\tau(\mathrm{tr}(X)\cdot
E_{1,1})=\mathrm{tr}(X)\tau(E_{1,1})$.


The following result by Dolinar et al. should be mentioned here:

\bigskip

\noindent\textbf{Theorem.}\textit{ (see \cite{Do}) Let }$K$\textit{ be a
field, and let }$\psi:\mathrm{M}_{n}(K)\rightarrow\mathrm{M}_{n}(K)$\textit{
be a bijective map which preserves the commutator Lie product. Then there is
an invertible matrix }$T\in\mathrm{M}_{n}(K)$\textit{, a field authomorphism
}$f$\textit{ of }$K$\textit{, and a function }$\tau:\mathrm{M}_{n}%
(K)\rightarrow K$\textit{, where }$\tau(X)=0$\textit{ for all matrices of
trace zero such that:}

\noindent\textit{(i) for }$n\geq3$\textit{ and }$K$\textit{ with a least
}$2^{n-1}$\textit{ elements, either}
\[
\psi(X)=TX^{f}T^{-1}+\tau(X)I\text{\textit{ for all }}X\in\mathrm{M}_{n}(K),
\]
\textit{or}%
\[
\psi(X)=-T(X^{f})^{\top}T^{-1}+\tau(X)I\text{\textit{ for all }}X\in
\mathrm{M}_{n}(K);
\]
\noindent\textit{(ii) for }$n=2$\textit{ and }$\mathrm{char}K\neq2$\textit{,}
\[
\psi(X)=TX^{f}T^{-1}+\tau(X)I\text{\textit{ for all }}X\in\mathrm{M}_{n}(K).
\]

\bigskip

Considering this theorem, we note that if we consider the functions
\[
\sigma_{1},\sigma_{2}: \mathrm{M}_{n}(K) \to\mathrm{M}_{n}(K),
\]
defined by
\[
\sigma_{1}(X) = TX^{f}T^{-1} \qquad\mathrm{and} \qquad\sigma_{2}(X) =
-T(X^{f})^{\top}T^{-1}
\]
for all $X\in\mathrm{M}_{n}(K)$, then by the foregoing constructions and
considerations, $\sigma_{1}$ is an automorphims of $\mathrm{M}_{n}(K)$ and
$\sigma_{2}$ is the negative of an anti-automorphism of $M_{n}(K)$ (in both
cases as rings), and as before, we have exact descriptions of them in the
terms of generators of $\mathrm{M}_{n}(K)$. However, we know nothing more
about $\tau$.



\bigskip

\noindent\textbf{Funding}

\bigskip

\noindent The second author was partially supported by the National Research,
Development and Innovation Office of Hungary (NKFIH) K119934. The research of
the fourth author was funded by the Polish National Science Centre Grant DEC-2017/25/B/ST1/00384.

\bigskip


\end{document}